\documentclass{amsart}

\usepackage{amsfonts, amssymb}

\usepackage{amssymb,amsmath,amscd,euscript,verbatim,array}
\usepackage{geometry}
\usepackage{anysize}
\usepackage{fancyhdr}
\usepackage{indentfirst}
\usepackage{graphicx}
\usepackage{color}
\usepackage{ifpdf}
\usepackage{marginnote}

\newcommand{\N}{\mathbb{N}}
\newcommand{\Z}{\mathbb{Z}}
\newcommand{\R}{\mathbb{R}}

\newcommand{\til}{\tilde}
\newcommand{\pa}{\partial}

\newcounter{item}[section]

\newcommand{\Proof}{\textbf{Proof}\hspace{0.3cm}}
\newcommand{\End}{\ensuremath{\hfill{\Box}}\\}

\newtheorem{Theorem}{Theorem}
\newtheorem{Definition}{Definition}[section]
\newtheorem{Proposition}[Definition]{Proposition}

\newtheorem{Lemma}[Definition]{Lemma}

\theoremstyle{definition}
\newtheorem{Remark}[Definition]{Remark}

\renewcommand{\th}{\theta}

\newcommand{\eps}{\varepsilon}

\newcommand{\T}{\mathbb{T}}

\begin{document}

\title[On the Destruction of Invariant Lagrangian Graphs for C.S. Twist Maps]{On the Destruction of Invariant Lagrangian Graphs for Conformal Symplectic Twist Maps}

\author{Alfonso Sorrentino}
\address{Department of Mathematics, University of Rome Tor Vergata, Via della ricerca scientifica 1, 00133 Rome, Italy}
\email{sorrentino@mat.uniroma2.it}

\author{Lin Wang}
\address{School of Mathematics and Statistics, Beijing Institute of Technology, Beijing 100081, China}
\email{lwang@bit.edu.cn}



\begin{abstract} {
In this article we investigate the fragility of invariant Lagrangian graphs for dissipative maps, focusing on their destruction under small perturbations.
Inspired by Herman's work on conservative systems,  we prove that all $C^0$-invariant Lagrangian graphs for an integrable dissipative  twist maps can be destroyed by perturbations that are arbitrarily small  in the $C^{1-\eps}$-topology. This result is sharp, as evidenced by the persistence of $C^1$-invariant graphs under $C^1$-perturbations guaranteed by the normally hyperbolic invariant manifold theorem.
}
\end{abstract}

\maketitle

\section{Introduction and main results}
Dissipative maps and flows have been a focal point of research in dynamical systems, yielding several profound results that have significantly advanced our understanding of their structure and behavior. One of the earliest milestones in this field was achieved by Birkhoff in 1932, who established the existence of attractors for dissipative twist maps~\cite{B}. Building on this foundation, subsequent works revealed the presence of rich and intricate dynamics, including periodic and quasi-periodic motions, KAM tori, Aubry-Mather sets, and more (see, for example,~\cite{CCD1, CCD2, Ca, L1,Ma,MS} and references therein). Further insights into the fine properties of the Birkhoff attractor were provided by Le Calvez~\cite{L2}, who explored its detailed characteristics and dynamical implications.
From a physical perspective, the study of dissipative twist maps has also been enriched by specific models illustrating the emergence of strange attractors, as highlighted in~\cite{Za}.  \
{Very recently, the notion of Birkhoff attractor has been extended to higher dimentional dissipative maps  in \cite{AHV}, and their properties have been thoroughly investigated by means of symplectic topological tools. We also refer the readers to \cite{AA} and \cite{AF} for a detailed investigation of the dynamical properties of  conformal symplectic systems and their invariant submanifolds. }\\

In this article, we focus on invariant Lagrangian graphs for dissipative maps, with particular attention to their fragility—namely, the possibility of their destruction by sufficiently small perturbations in a suitable topology.\\

Let us begin by describing the conservative case, where the destruction of invariant Lagrangian graphs for twist maps and Hamiltonian systems has been extensively studied, being viewed as a natural counterpart to KAM theory.

A seminal contribution by Herman~\cite{H1} demonstrated that any sufficiently small $C^{3-\varepsilon}$ perturbation of a twist map of the annulus can destroy an invariant circle with a prescribed frequency, highlighting the sharpness of the regularity assumptions in KAM theory~\cite{H33}. Advancing this line of research, Mather~\cite{M4} showed that for any Liouville frequency, an invariant circle can be destroyed by an arbitrarily small $C^\infty$ perturbation. Forni~\cite{F} further proved that even in the real-analytic ($C^\omega$) topology, an invariant circle with a frequency belonging to a proper subset of the Brjuno numbers can be eliminated under small perturbations. Notably, Takens~\cite{T} provided the first example where all invariant circles are annihilated by a $C^1$-small perturbation.

Turning to the destruction of invariant Lagrangian tori for symplectic twist maps of $\mathbb{T}^d \times \mathbb{R}^d$, Herman~\cite{H3} showed that any small $C^{d+1-\varepsilon}$ perturbation suffices to destroy all such tori, with the same conclusion holding for real-analytic systems~\cite{W2}. An alternative approach, developed in~\cite{SW}, leverages Percival’s Lagrangian to achieve the destruction of all invariant tori in the context of the generalized Frenkel-Kontorova model. \\

The difference between the conservative and dissipative cases is significant in many perspectives.  For example, in the dissipative setting, at most one invariant graph may exist, which can be interpreted as a particular type of Birkhoff attractor.

When considering nearly integrable systems, the fundamental difference between conservative and dissipative systems lies in the fact that the invariant tori of dissipative integrable systems (see (\ref{F1}) and (\ref{fbeta}) below for instance) are normally hyperbolic. The normally hyperbolic invariant manifold theorem \cite{HPS} guarantees the persistence of $C^1$-invariant graphs under $C^1$-perturbations, which leads to significant differences in the preservation and breakdown of their invariant tori under small perturbations.

The existence of invariant graphs in conservative nearly integrable systems is guaranteed by the KAM theorem, where the smoothness of the perturbation is related to the frequency of the invariant graphs and the degrees of freedom of the system. To the best of our knowledge, there remains a gap between the results on the existence and non-existence of invariant graphs. Specifically, for \( d = 1 \), small perturbations in the \( C^3 \) topology \cite{H33} ensure the persistence of invariant circles, whereas perturbations leading to the breakdown of all invariant circles require smallness in the \( C^{2-\varepsilon} \) topology \cite{H1}---whether this can be further improved remains unclear. For \( d \geq 2 \), the gap in related results is even larger. For Hamiltonian functions, small perturbations in the \( C^{2d+w} \) topology (where \( w \) is a certain modulus of continuity \cite{TL}) ensure the persistence of invariant graphs, while small perturbations in the \( C^{d+1-\varepsilon} \) topology lead to the breakdown of all invariant graphs.

In this paper, we prove that at all of $C^0$-invariant Lagrangian graphs for a family of integrable conformal symplectic twist maps can be destroyed
by an arbitrarily small perturbation in the $C^{1-\eps}$-topology. As a consequence, this provides an evidence the the Birkhoff attractor can fail to be a graph under an arbitrarily small perturbation in the $C^{1-\eps}$ topology. This result is sharp in view of the normally hyperbolic invariant manifold theorem. Furthermore,  it was  pointed out in \cite[Remark 2, Page 52]{HPS} that a Lipschitz invariant graph persists under sufficiently small Lipschitz perturbations by  using tangent cones instead of tangent planes.
\\

Let us now describe our main results in more details.\\

\subsection{$1$-dimensional case}

Let us first start with the $1$-degree of freedom case; we denote by $\T:=\R/\Z$. \\

\begin{Definition}[Dissipative twist map of the cylinder]\label{dissmap}
A {\it dissipative twist map}  of the (infinite) cylinder is a  $C^1$ diffeomorphism $f:\T \times \R \rightarrow \T \times \R  $ that admits a lift
$F:\R^{2}\rightarrow \R^2$, $F(x,y)=(X(x,y),Y(x,y))$, satisfying the following conditions:
\begin{itemize}
    \item[(i)]  $F$ is isotopic to the identity;
    \item[(ii)] ({\it Twist condition}) The map $\psi:(x,y)\mapsto(x,X(x,y))$ is a diffeomorphism of $\R^2$;
    \item[(iii)]    ({\it Dissipative condition}) There exists $\lambda\in (0,1)$ such that for all $(x,y)\in \R^2$
\[0<\delta\leq \det(DF(x,y))\leq \lambda.\]
\end{itemize}
\end{Definition}

\medskip

\begin{Remark} \label{remark1}
{\bf (i)} {Observe that $F(x +1, y)=F(x, y)+(1, 0)$ for every $(x,y)\in \R^2$}. Clearly, any $C^1$ map $F: \R^2\rightarrow \R^2$ satisfying (i)-(iii) induces, by projection, a dissipative map  $f$: $\T\times\R\mapsto \T\times\R$. In the following, we will equivalently consider the lift $F$ or the map $f$.\\
{\bf (ii)} One could also consider maps defined on the finite cylinder (or annulus) $\T\times [a,b]$, with $a<b$, up to suitably adapt the definition of dissipative map to this case. Although our result  clearly extends to this case, for the sake of simplicity we only discuss the case of the infinite cylinder.\\
{\bf (iii)} We remark that the dissipative condition (iii) in  Definition \ref{dissmap} is very general and, for example, does not require the determinant of the Jacobian of the map to be constant (differently from what happens in the higher dimensional case, see the definition of conformal symplectic map in Definition \ref{confsympmap}).\\
\end{Remark}

In this note, we investigate the stability of invariant graphs.

\begin{Definition}[$C^0$-invariant graph]\label{c0graph}
$\Gamma\subset\T\times\R$ is called a $C^0$-invariant graph of $f$ if
\begin{itemize}
\item[(i)] $\Gamma=\{(\theta,\psi(\theta)):\; \theta \in \T\}$, where $\psi:\T\to\R$ is a continuous function;
\item[(ii)] $\Gamma$ is invariant under the action of  $f$.
\end{itemize}
\end{Definition}

\begin{Remark}
{(i)} By the twist condition, if $f$ is of class $C^1$, then $\psi$ is a Lipschitz function on $\T$ (see \cite[Proposition 2.2]{H1}).\\
{(ii)} Equivalently, if $F: \R^2 \to \R^2$ denotes a lift of $f$ and $\tilde{\psi}:\R\to\R$  a $1$ lift of $\psi$ (which is a $1$-periodic function on $\R$), then  the graph $\tilde{\Gamma}:=\{(x,\tilde{\psi}(x)) |\ x\in\R\}$  is invariant by $F$.
\end{Remark}

In the following, it will be more convenient to consider the lifts of the maps and the graphs.\\

Let $\alpha:=(\alpha_1,\alpha_2)\in \R^2$, $\lambda\in (0,1)$. We consider $F_{\lambda,\alpha}$:
\begin{equation}\label{F1}
F_{\lambda,\alpha}(x,y):=(x+\alpha_1+\lambda y,\alpha_2+\lambda y).
\end{equation}
The parameter $\lambda$ controls the dissipation, and $\alpha_1,\alpha_2$
 are constants that determine the translation in the $x$ and $y$ directions, respectively.
It is easy to check that this map induces a dissipative twist map. In particular, one can completely describe its dynamics:  the graph $\Gamma_{\lambda,\alpha}:= \R \times \{\frac{\alpha_2}{1-\lambda}\}$ is invariant under the action of $F_{\lambda,\alpha}$ (the dynamics on it reduces to a translation), while
any other orbit of $F_{\lambda,\alpha}$ is forward-asymptotic to $\Gamma_{\lambda,\alpha}$. In fact, if we denote
if $(x_k,y_k):= F^k_{\lambda,\alpha}(x_0,y_0)$ for $k\geq 0$, then
$$
y_k = \alpha_2 \sum_{j=0}^k \lambda^j + \lambda^k y_0 \longrightarrow \frac{\alpha_2}{1-\lambda} \qquad {\rm as\;} k\to +\infty.
$$

\medskip
We want to consider a perturbation of this family of maps by a $C^\omega$  function $\varphi:\R\to\R$ which is $1$-periodic and satisfies $\int_{0}^1\varphi(x)dx=0$, namely
\[F^{\flat}_{\lambda,\alpha}(x,y)=(x+\alpha_1+\lambda y+\varphi(x),\alpha_2+\lambda y+\varphi(x)).\]
Clearly, $\det DF_{\lambda,\alpha}^\flat(x,y)=\lambda$ for all $(x,y)\in \R^2$, hence this induces a family of dissipative twist maps.\\

\medskip

\bigskip

We can now state our first result for dissipative twist maps.

\begin{Theorem}\label{Mt}
Given  $\lambda\in (0,1)$ and $0<\eps\ll 1$, there exists a sequence of  trigonometric polynomials $\{\varphi^{\lambda}_n\}_{n\in \N}$ such that all of the $C^0$-invariant graphs for the family of $\{F_{\lambda,\alpha}\}_{\alpha\in\R^2}$ can be destroyed
by perturbing the maps with $\{\varphi^{\lambda}_n\}_{n\in \N}$. Moreover,
 \begin{itemize}
 \item [(I)]the degree $N$ of $\varphi^{\lambda}_n$ satisfies $N=O(n^{1+\eps})$ as $n\to \infty$;\\
 \item [(II)]  $\|\varphi^{\lambda}_n\|_{C^{1-\eps}}=O(\frac{1}{n^\eps})$  as $n\to \infty$;\\
 \item [(III)] for every $n\in \N$,  we have $\|\varphi^{\lambda}_n\|_{C^1}\leq 1$ for all $\lambda\in (0,1)$, and $\|\varphi^{\lambda}_n\|_{C^1}=O(1-\lambda)$ as $\lambda\to 1^-$.
 \end{itemize}
\end{Theorem}

Theorem \ref{Mt} shows the fragility of $C^0$-invariant Lagrangian graphs, and it also implies the Birkhoff attractor might fail to be a graph up to an arbitrarily small perturbation in the
$C^{1-\eps}$-topology. \\

\subsection{$d$-dimensional case}
We want now to describe a higher dimensional analogue of Theorem \ref{Mt}. Let us introduce the right setting.
Denote $\T^d:=\R^d/\Z^d$.

\begin{Definition}[conformal symplectic twist map]\label{confsympmap}
A {\it conformal symplectic twist map}  of the $d$-dimensional cylinder is a  $C^1$ diffeomorphism $f:\T^d \times \R^d \rightarrow \T^d \times \R^d  $ that admits a lift
$F:\R^{2d}\rightarrow \R^{2d}$, $F(x,y)=(X(x,y),Y(x,y))$, satisfying:
\begin{itemize}
    \item[(i)] $F$ is isotopic to the identity;
      \item[(ii)] ({\it Twist condition}) The map $\psi:(x,y)\mapsto(x,X(x,y))$ is a diffeomorphism of $\R^{2d}$;
    \item[(iii)]  ({\it conformal symplectic condition}) There exists $\lambda\in (0,1)$ such that
    $$
    F^*(dx\wedge dy) = \lambda\, dx\wedge dy
    $$
    where $F^*$ denotes the pull-back by $R$; namely the map $F$ rescales the canonical symplectic form by a constant factor. Equivalently, this could be rephrased by saying that
    the 1-form $\alpha:=F^*ydx-\lambda ydx$ is a closed $1$-form.
    In particular, if there exists $S:\R^{2d}\to\R$ such that $\alpha=dS$, then $F$ is called exact, and the function $S$ is called a generating function of $F$. Namely, $F$
is generated by the following equations
\begin{equation*}
\begin{cases}
\lambda y=-\partial_1 S(x,X),\\
Y=\partial_2 S(x,X),
\end{cases}
\end{equation*}
where $F(x,y)=:(X,Y)$.

\end{itemize}

\end{Definition}

\begin{Remark}
{\bf (i)} Similarly to what remarked in Remark \ref{remark1} (i), {one has that $F(x +m, y)=F(x, y)+(m, 0)$ for every $m\in \Z^d$ and $(x,y)\in \R^{2d}$}. Moreover,
the projection of any $C^1$ map $F: \R^{2d}\rightarrow \R^{2d}$ satisfying (i)-(iii) in Definition \ref{confsympmap}, induces, by projection, a conformal symplectic twist map $f$: $\T^d\times\R^d\mapsto \T^d\times\R^d$. In the following, we will equivalently consider the map $f$ or its lift $F$.\\
{\bf (ii)} For {$d\geq 2$}, the definition of  conformal symplectic twist maps is more restrictive than the one we gave in Definition \ref{dissmap}. In fact, Definition \ref{confsympmap} implies that  $\mathrm{det}(DF(x,y))$ must be a constant (see \cite{Li}).
\end{Remark}

Also in this setting we want to investigate the stability of invariant graphs. However,  in dimension $d\geq 2$   invariant graphs that are more relevant from a dynamical system point of view, are those satisfying the property fo being {\it Lagrangian}  (this additional property is automatically satisfied in dimension $d=1$).\\

\begin{Definition}[$C^0$-invariant Lagrangian graph]\label{lag}
$\mathcal L \subset\T^d\times\R^d$ is called a $C^0$-invariant Lagrangian graph of $f$ if
\begin{itemize}
\item[(i)] $\mathcal L =\{(\theta,\psi(\theta)):\; \theta \in \T\}$, where $\psi=(\psi_1,\ldots, \psi_d):\T^d\to\R^d$ is a continuous function;
\item[(ii)] $\mathcal L$ is invariant under the action of  $f$;
\item[(iii)] the $1$-form $\Psi(\theta)\;d\theta:=\sum_{i=1}^d\psi_i(\theta)d\theta_i$, where $\th:=(\th_1,\ldots,\th_d)\in \T^d$, is closed in the sense of distribution.
Equivalently, there exists a $C^1$ function $\eta:\T^d\to\R$ such that $\Psi(\theta)=c+D\eta(\th)$ for all $\theta\in \T^d$, where $c\in \R^d$ is a constant vector and $D\eta:=(\frac{\partial \eta}{\partial \th_1},\dots,\frac{\partial \eta}{\partial \th_d})$.\\
\end{itemize}
\end{Definition}

\begin{Remark}
{\bf (i)} For $d=1$, Definition \ref{lag} simplifies Definition \ref{c0graph}, since a homotopically nontrivial invariant curve is automatically Lagrangian  in this case. \\
{\bf (ii)} If $\mathcal L$ is a $C^1$ graph, this definition reduces to the classical one, namely every tangent space is a Lagrangian subspace (see, for example, \cite{CdS}). In particular, it is a well-known fact that $C^1$ Lagrangian graphs in the cotangent bundle of a manifold correspond to the graph of closed $1$-forms (see for example \cite[Section 3.2]{CdS}).
\end{Remark}

\medskip

Let $\beta\in \R^d$ and $\lambda \in (0,1)$. Consider the map $F_{\lambda,\beta}:\R^d\times\R^d\to\R^d\times\R^d$ given by
\begin{equation}\label{fbeta}
F_{\lambda,\beta}(x,y)=(x+\lambda\beta+\lambda  y,\lambda y).
\end{equation}

This map clearly induces a conformal symplectic twist map of $\T^d\times \R^d$. In particular, $\mathcal L_0 := \R^d\times \{0\}$
is invariant under the action of $F_{\lambda,\beta}$ (the dynamics on it reduces to a translation by $\lambda \beta$ ), while
any other orbit of $F_{\lambda,\beta}$ is forward-asymptotic to $\mathcal L_0$. In fact, if we denote
if $(x_k,y_k):= F^k_{\lambda,\beta}(x_0,y_0)$ for $k\geq 0$, then
$$
y_k =  \lambda^k y_0 \longrightarrow 0 \qquad {\rm as\;} k\to +\infty.
$$

We want to consider a perturbation of this family of maps by a $C^\omega$  function $\Phi:\R^d\to\R$ which is $\Z^d$-periodic, namely
\[F_{\lambda,\beta}^\flat(x,y)=(x+\lambda (\beta+y+D\Phi(x)),\lambda (y+D\Phi(x))).\]
One can easily check that $F_{\lambda,\beta}^\flat$ are still conformal symplectic (see section \ref{subsechermantd} for more details).\\

\begin{Remark}
The perturbation of the map $F_{\lambda,\beta}$ can be more naturally described in terms of a perturbation of its generating function (see Definition \ref{confsympmap} (iii)). More specifically,
$F_{\lambda,\beta}$ is generated by
$$S_{\lambda,\beta}(x,X):=\frac{1}{2}\langle X-x,X-x\rangle-\lambda \langle \beta,X-x\rangle$$
where $\langle \cdot,\cdot\rangle$ denotes the standard inner product in $\R^d$. One can easily check that $F_{\lambda,\beta}^\flat$ is generated by
$$
S_{\lambda,\beta}^\flat(x,X):=\frac{1}{2}\langle X-x,X-x\rangle-\lambda \langle \beta,X-x\rangle+\lambda \Phi(x).
$$
Since $S_{\lambda,\beta}^\flat(x,X)$ is defined up to additive constants, we can assume that $\Phi$ has zero average.
\end{Remark}

We can prove the following result, which is an analogue of Theorem \ref{Mt} in the higher dimensional case.
\begin{Theorem}\label{Mt2}
Given  $\lambda\in (0,1)$ and $0<\eps\ll 1$,  there exists a sequence of  trigonometric polynomials $\{\Phi^{\lambda}_n\}_{n\in \N}$ such that all of the $C^0$-invariant graphs for the family  of $\{S_{\lambda,\beta}\}_{\beta\in\R^d}$ can be destroyed
by perturbing the maps with $\{\Phi^{\lambda}_n\}_{n\in \N}$. Moreover,
 \begin{itemize}
 \item [(I)] the degree $N$ of $\Phi^{\lambda}_n$ satisfies $N=O(n^{\frac{1}{d}+\eps})$ as $n\to \infty$;\\
 \item [(II)]  $\|\Phi^{\lambda}_n\|_{C^{2-\eps}}=O(\frac{1}{n^\eps})$  as $n\to \infty$;\\
 \item [(III)]  for every $n\in \N$,  we have $\|\Phi^{\lambda}_n\|_{C^2}\leq d$ for all $\lambda\in (0,1)$,  and $\|\Phi^{\lambda}_n\|_{C^2}=O(1-\lambda)$ as $\lambda\to 1^-$.
 \end{itemize}
\end{Theorem}

\medskip

\begin{Remark}\label{RemarkmatrixA1}
Let $A$ be a $d\times d$ symmetric positive definite  matrix, and let $\beta\in \R^d$ be a constant vector. By a similar argument, one can show that Theorem \ref{Mt2} also holds for the family of maps $F_{\lambda,\beta}$ generated by
\[S_{\lambda,\beta}(x,X):=\frac{1}{2}\langle X-x,A(X-x)\rangle-\lambda \langle \beta,A(X-x)\rangle.\]
In fact, one can construct the perturbed maps $F_{\lambda,\beta}^\flat$ generated by
\[S_{\lambda,\beta}(x,X):=\frac{1}{2}\langle X-x,A(X-x)\rangle-\lambda \langle \beta,A(X-x)\rangle+\lambda W(x),\]
where $W:\R^d\to \R$ is a  $\Z^d$-periodic function of class $C^\omega$ satisfying the following:
\begin{itemize}
\item $\int_{\T^d}W(x)dx=0$;
\item there exists a  $\Z^d$-periodic $\Phi:\R^d\to \R$ of class $C^\omega$ such that  $DW(x)=AD\Phi(x)$ for each $x\in \R^d$.
\end{itemize}
Accordingly, $W(x)$ is determined by the Fourier coefficients of $\Phi$ and we also have $\int_{\T^d}\Phi(x)dx=0$. Then
\[F_{\lambda,\beta}^\flat(x,y)=\left(x+\lambda (\beta+D\Phi(x)+A^{-1}y),\lambda (y+AD\Phi(x))\right).\]
{See also Remark \ref{RemarkmatrixA2} for more indications on how to modify the argument yielding the proof of an analog of Theorem \ref{Mt2} in this setting.}\\

\end{Remark}

\bigskip

\subsection*{Concluding remarks} \phantom{x}

\smallskip

\noindent {\bf (i)} In general, dissipative twist maps may admit invariant curves that fail to be graphs (see~\cite[Proposition 15.3]{L2}).
However, the perturbation constructed in~\cite{L2} is not small in the $C^0$-topology.
Consequently, it remains an open question whether such non-graph invariant curves can persist in the perturbed systems described in Theorems~\ref{Mt} and~\ref{Mt2}.

\noindent {\bf (ii)} It is straightforward to verify that the set $\mathbb{T} \times \left\{\frac{\alpha_2}{1-\lambda}\right\}$ (respectively, $\mathbb{T}^d \times \{0\}$) is normally hyperbolic under the map $F_{\lambda,\alpha}$ (resp.\ $F_{\lambda,\beta}$).
By the theory of normally hyperbolic invariant manifolds \cite{HPS}, a $C^1$-smooth invariant graph persists under sufficiently small $C^1$-perturbations.
This demonstrates that Theorems~\ref{Mt} and~\ref{Mt2} are sharp with respect to the topology of the perturbation.
This is different from what happens in the conservative case considered in \cite{H3}, where the unperturbed map has a foliation by invariant tori (hence, they are not normally hyperbolic) and  the perturbation that destroys them can be taken small in the $C^{d+1-\varepsilon}$-topology.

For perturbations with specific structure, such as the dissipative standard map
\[
\varphi(x) = \frac{k}{2\pi} \sin(2\pi x),
\]
it was shown in \cite{Bo} that invariant graphs can be destroyed when $k > k_0:=\frac{2(1+\lambda)}{2+\lambda}$, where the value of $k_0$ can be also obtained by letting $-m=M$ in (\ref{mM2}) below.
This generalizes (by sending $\lambda\to 1^-$) a celebrated result by Mather \cite{M1} for the conservative case, where the critical threshold is $k > \frac{4}{3}$.

\smallskip
\noindent {\bf (iii)} The result presented in \cite{CCD1} concerns a family of conformal symplectic maps \( f_\mu \) on the torus \( \mathbb{T}^d \) for \( d \geq 2 \). The map \( f_\mu \) is exact if and only if \( \mu = 0 \). Using an {\it a posteriori} format, the authors prove that if \( f_{\mu_0} \) admits an invariant torus, then there exists another parameter \( \mu_e \) such that the torus persists for \( f_{\mu_e} \). Here, the parameter \( \mu \) is referred to as the {\it drift parameter}. In our \( d \)-dimensional setting, the requirement of exactness implies that \( \mu \) must be necessarily zero.

\medskip

\subsection*{Organization of the paper}
{ The paper is organized as follows.  Section \ref{S1} is focused on the $1$-dimensional setting and it presents the proof of Theorem \ref{Mt}. More specifically, in Subsection \ref{Herman1} we prove one of the main tools our analysis, namely an analogue Herman's {\it a-posteriori equality} for the existence of invariant Lagrangian graph (see Proposition \ref{hf1}), extending to the dissipative setting the result from \cite{H1}. This will be exploited in Subsection \ref{subsect2.2} to construct the pertubation, using tools of approximation theory ({\it i.e.}, Jackson's approximation method), and complete the proof of Theorem \ref{Mt}.  In Section \ref{S2}, inspired by the foundational contributions of Herman \cite{H2,H3}, we extend the previous analysis to the higher dimensional case; in Subsection \ref{subsechermantd} we discuss Herman's formula for higher dimensional conformal symplectic maps, while in Subsection \ref{subsechermantd} we construct the perturbation and complete the proof of Theorem \ref{Mt2}. Finally, in Appendix \ref{appendix}, for the reader's sake, we present a proof of Jackson's approxation theorem using an argument similar to  \cite[Theorem 2.12]{A}.\\
}

\medskip

\noindent\textbf{Acknowledgement.} The authors would like to thank Patrice Le Calvez for his valuable comments and Qinbo Chen for pointing out a mistake in a preliminary version.  LW is supported by NSFC Grant No. 12122109.
AS acknowledges the support of the Italian Ministry of University and Research's PRIN 2022 grant ``Stability in Hamiltonian dynamics and beyond'', as well as the Department of Excellence grant MatMod@TOV (2023-27) awarded to the Department of Mathematics of University of Rome Tor Vergata. AS is a member of the INdAM research group GNAMPA and the UMI group DinAmicI.

\medskip

\section{Proof of Theorem \ref{Mt}}\label{S1}
We recall that $F^{\flat}_{\lambda, \alpha}$ denotes the perturbed map
\[F^{\flat}_{\lambda, \alpha}(x,y):=(x+\alpha_1+\lambda y+\varphi(x),\alpha_2+\lambda y+\varphi(x)),\]
where $\varphi:\R\to\R$ is a $1$-periodic $C^1$  function satisfying $\int_{0}^1\varphi(x)dx=0$. It is clear that det$DF_{\lambda, \alpha}^\flat(x,y)=\lambda$ for all $(x,y)\in \R^2$.

\subsection{Herman's formula for dissipative twist maps}\label{Herman1}
In this section we are going to prove an analogue of Herman's formula \cite{H1} in this dissipative setting. More specifically, this consists in an  {\it a-posteriori equality} that must be satisfied if there exists invariant Lagrangian graph.

\begin{Proposition}\label{hf1}
$F^\flat_{\lambda,\alpha}$ admits a $C^0$-invariant graph $\tilde{\Gamma}:=\{(x,\tilde{\psi}(x))\ |\ x\in \R\}$ if and only if
\begin{equation}\tag{A}\label{a}
\frac{1}{1+\lambda}g(x)+\frac{\lambda}{1+\lambda}g^{-1}(x)=x+\frac{1}{1+\lambda}((1-\lambda)\alpha_1+\lambda\alpha_2+\varphi(x)) \qquad \forall\; x\in \R
\end{equation}
where $g(x):=x+\alpha_1+\lambda\tilde{\psi}(x)+\varphi(x)$.
\end{Proposition}
\begin{Remark}
For $\lambda=1$, the formula (A) was established by Herman in 1983 (see \cite[Section 2.4]{H1}).
\end{Remark}

\Proof
For simplifying the notation, we denote $F:=F^\flat_{\lambda,\alpha}$. \\
$(\Longrightarrow)$ Let us assume that  $F$ admits an invariant graph $\tilde{\Gamma}:=\{(x,\tilde{\psi}(x))\ |\ x\in \R\}$. Then
\begin{align*}
F(x,\tilde{\psi}(x))&=(x+\alpha_1+\lambda\tilde{\psi}(x)+\varphi(x),\alpha_2+\lambda\tilde{\psi}(x)+\varphi(x))\\
&=(g(x),\tilde{\psi}(g(x))),
\end{align*}
where
\begin{equation}\label{g}
g(x)= x +\alpha_1+\lambda \tilde{\psi}+\varphi.
\end{equation}
Note that
\begin{align*}
F^{-1}(x,\tilde{\psi}(x))&=(x-\tilde{\psi}(x)-\alpha_1+\alpha_2,\frac{1}{\lambda}(\tilde{\psi}(x)-\alpha_2-\varphi(x-\tilde{\psi}(x)-\alpha_1+\alpha_2))\\
&=(g^{-1}(x),\tilde{\psi}(g^{-1}(x))),
\end{align*}
which implies
\[g^{-1}(x)= x -\tilde{\psi}(x)-\alpha_1+\alpha_2.\]
Combining with (\ref{g}), we have formula (\ref{a}).\\

\noindent $(\Longleftarrow)$ Conversely, let us assume that (\ref{a}) holds. We only need to show
\begin{equation}\label{c}
\tilde{\psi}\circ g=\alpha_2+\lambda \tilde{\psi}+\varphi.
\end{equation}
By (\ref{a}),
\begin{equation}\label{b}
g\circ g+\lambda\,\mathrm{Id}=(1+\lambda)g+(1-\lambda)\alpha_1+\lambda\alpha_2+\varphi\circ g.
\end{equation}
From the definition of $g$, we have
\[g\circ g=g+\alpha_1+\lambda\tilde{\psi}\circ g+\varphi\circ g.\]
Combining with (\ref{b}), we obtain
\[\tilde{\psi}\circ g=g-\alpha_1+\alpha_2-\mathrm{Id},\]
By using the definition of $g$ again, we obtain (\ref{c}).
\End

Since $F$ is of class $C^1$, $\tilde{\psi}$ is Lipschitz by \cite[Proposition 2.2]{H1}. Then $g$ is also Lipschitz. Denote
\[G:=\max\{\|Dg\|_{L^\infty},\|Dg^{-1}\|_{L^\infty}\}.\]
Let $\mathfrak{D}$ be the set of  points of differentiability of $g$,
then $\mathfrak{D}$ has full Lebesgue measure on $\R$. By the construction of $F^\flat_{\lambda,\alpha}$, for each $x\in \mathfrak{D}$, $Dg(x),\ Dg^{-1}(x)>0$. We differentiate
(\ref{a}) and obtain
\begin{equation}\label{deri}
\frac{1}{1+\lambda}Dg(x)+\frac{\lambda}{1+\lambda}(Dg)^{-1}(g^{-1}(x))=1+\frac{1}{1+\lambda}D\varphi(x).
\end{equation}
Let
\[M:=\max_{x\in [0,1]}D\varphi(x),\quad m:=\min_{x\in [0,1]}D\varphi(x).\]
Since $\varphi$ is $1$-periodic, then $\int_0^1D\varphi dx=0$. Therefore, $m\leq 0\leq M$. Note that for $x\in \mathfrak{D}$
\[(Dg)^{-1}(g^{-1}(x))\geq \frac{1}{G},\quad Dg(x)=\frac{1}{Dg^{-1}(g(x))}\geq \frac{1}{G}.\]
By (\ref{deri}), we have
\begin{equation}\label{ggGG}
\frac{1}{G}\leq 1+\frac{1}{1+\lambda}m,
\end{equation}
which implies $m>-1-\lambda$.
The following part is divided into two cases.
\begin{itemize}
\item {\bf{Case 1:}} $G=\|Dg\|_{L^\infty}$;
\item {\bf{Case 2:}} $G=\|Dg^{-1}\|_{L^\infty}$.
\end{itemize}

\subsubsection{Case 1}

By (\ref{deri})
\begin{equation}\label{G}
\frac{1}{1+\lambda}G+\frac{\lambda}{1+\lambda}\frac{1}{G}\leq 1+\frac{1}{1+\lambda}M.
\end{equation}
A direct calculation shows
\begin{equation}\label{GM}
G\leq \frac{1+\lambda}{2}+\frac{M}{2}+\left(\frac{(1-\lambda)^2}{4}+\frac{M^2}{4}+\frac{(1+\lambda)M}{2}\right)^{\frac{1}{2}},
\end{equation}
which together with (\ref{ggGG}) implies
\begin{equation}\label{mM}
\frac{1}{1+\frac{1}{1+\lambda}m}\leq \frac{1+\lambda}{2}+\frac{M}{2}+\left(\frac{(1-\lambda)^2}{4}+\frac{M^2}{4}+\frac{(1+\lambda)M}{2}\right)^{\frac{1}{2}}.
\end{equation}
If $M\to 0^+$, then for $m\in (-1-\lambda,0)$,
\begin{equation}\label{mmMM}
-m\leq \frac{1+\lambda}{1-\lambda}M+O(M^2).
\end{equation}
In fact,  for the left hand side of (\ref{mM}), we have, by Taylor's formula,
\[\frac{1}{1+\frac{1}{1+\lambda}m}=1-\frac{1}{1+\lambda}m+\frac{1+\lambda}{(1+\lambda+\xi)^3}m^2> 1-\frac{1}{1+\lambda}m,\]
where $\xi\in (m,0)$.
For the right hand side of (\ref{mM}), we have, as $M\to 0^+$
\[\left(\frac{(1-\lambda)^2}{4}+\frac{M^2}{4}+\frac{(1+\lambda)M}{2}\right)^{\frac{1}{2}}\leq \frac{1-\lambda}{2}+\frac{1+\lambda}{2(1-\lambda)}M+O(M^2).\]
This verifies (\ref{mmMM}).

\subsubsection{Case 2}
By (\ref{deri})
\begin{equation}\label{G2}
\frac{1}{1+\lambda}\frac{1}{G}+\frac{\lambda}{1+\lambda}G\leq 1+\frac{1}{1+\lambda}M.
\end{equation}
A direct calculation shows
\begin{equation}\label{GM2}
G\leq \frac{1}{2\lambda}\left(1+\lambda+{M}+\left((1+\lambda+{M})^2-4\lambda\right)^{\frac{1}{2}}\right),
\end{equation}
which together with (\ref{ggGG}) implies for $m\in (-1-\lambda,0)$
\begin{equation}\label{mM2}
\frac{1}{1+\frac{1}{1+\lambda}m}\leq \frac{1}{2\lambda}\left(1+\lambda+{M}+\left((1+\lambda+{M})^2-4\lambda\right)^{\frac{1}{2}}\right).
\end{equation}
We take $M=\frac{1}{n}$ and let $n\to +\infty$,
\begin{equation}\label{mmMM2}
-m\leq 1-\lambda^2+\frac{\lambda(1+\lambda)}{1-\lambda}{\frac{1}{n}}+O\left(\frac{1}{n^2}\right).
\end{equation}

\medskip
It is clear to see that if (\ref{mmMM2}) is invalid, then (\ref{mmMM}) does not hold. Therefore, in order  to destroy invariant graphs, it suffices to construct a perturbation satisfying the following {\bf criterion}
\begin{itemize}
\item [($\diamondsuit$)]   $-m <1+\lambda$ and $M\leq \frac{1}{n}$ such that (\ref{mmMM2}) is invalid as $n\to+\infty$.
\end{itemize}
\subsection{Construction of the perturbation}\label{subsect2.2}
We start this part by a lemma.
\begin{Lemma}\label{Csp}
Given   $0<\eps\ll 1$,  there exists a sequence of $\Z^d$-periodic trigonometric polynomials $\{\phi_n\}_{n\in \N}$ such that
 \begin{itemize}
 \item $\int_{\T^d}\phi_n(x)dx=0$ and \[-\min_{x\in \T^d}\phi_n(x)=\Delta>0,\quad 0<\max_{x\in \T^d}\phi_n(x)\leq \frac{1}{n},\]
 where $\Delta$ is a constant independent of $n$,
 \item the degree $N:=N(n)$ of $\phi_n$ satisfies $N=O(n^{\frac{1}{d}+\eps})$ and $\|\phi_n\|_{C^{1-\eps}}=O(n^{-\eps})$  as $n\to \infty$.
 \end{itemize}
\end{Lemma}

\Proof
To fix notations, we choose $[0,1]^d$ as the fundamental domain of $\T^d$. Denote $x:=\{x_1,\ldots,x_d\}$. First, we construct a sequence of $\Z^d$-periodic $C^\infty$ functions  denoted by $\{f_n\}_{n\in \N}$. The trigonometric polynomial perturbation  can be then obtained by using the Jackson approximation and tools from complex analysis {(see below and Appendix \ref{appendix})}.

\smallskip

 We {assert} the existence of $\{f_n\}_{n\in \N}$ satisfying
\begin{itemize}
\item [(1)] $\int_{[0,1]^d}f_n(x)dx=0$;
\item [(2)] \[-\min_{x\in [0,1]^d} f_n(x)= \Delta+\frac{1}{4n},\quad \max_{x\in [0,1]^d} f_n(x)=\frac{3}{4n}.\]
\end{itemize}


\vspace{1em}

Let us prove the existence of this family of functions $\{f_n\}_{n\in \N}$. For simplicity, we still use $a$ to denote the $d$-dimensional coordinate $(a,\ldots,a)$. Given $n\in\N$, we require $f_n$ to satisfy the following conditions:
\begin{itemize}
\item on the interval $[0,1]^d$,
\[\max f_n(x)=f_n\big(\frac{1}{4}\big)=\frac{1}{n}-\frac{1}{4n}=\frac{3}{4n},\]
\[ -\min f_n(x)=-f_n\big(\frac{3}{4}\big)=\Delta+\frac{1}{4n}; \]
\item $f_n(0)=f_n(\frac{1}{2})=f_n(1)=0$, $f_n>0$ on $(0,\frac{1}{2})^d$;
\item  $f_n\leq 0$ on $(\frac{1}{2},1]^d$, and $f_n$ is supported on the $d$-dimensional cube
$\left[\frac{3}{4}-b_n,\frac{3}{4}+b_n\right]^d$, where
\[b_n:=\frac{1}{4}\left(\Delta+\frac{1}{4n}\right)^{-\frac{1}{d}}\left(\frac{3}{4n}\right)^{\frac{1}{d}}.\]
\end{itemize}

\smallskip

\textit{Heuristic argument:}
If we only require \( f_n \) to be Lipschitz, then it can be constructed by forming the graph of \( f_n \) as two hyperparallelepipeds:
\begin{itemize}
    \item the left (upward) one has the  height \( \frac{3}{4n} \) and the base of volume \( \left(\frac{1}{2}\right)^d \).
    \item the right (downward) one  has the height  \( \Delta+\frac{1}{4n} \) and base of volume \( (2b_n)^d \).
\end{itemize}
A direct calculation shows that \( \int_0^1 f_n(x) \, dx = 0 \). We then modify \( f_n \) into a \( C^\infty \) function while ensuring that \( \int_0^1 f_n(x) \, dx = 0 \) remains valid.

\vspace{1em}

In order to use Jackson's approximation, we have to estimate $\|f_n\|_{C^k}$ as $n\to \infty$ for each $k\in \N$. Based on elementary Fourier analysis, we have
\[\|f_n^{(m)}\|_\infty\leq  C_{d,\lambda}\left(\frac{1}{4\pi b_n}\right)^m,\]
where $C_{d,\lambda}$ is a constant only depending on $d$ and $\lambda$.
It follows that
\begin{equation}\label{fck}
\|f_n\|_{C^k}\leq C_{d,\lambda}\sum_{m=0}^k\left(\frac{1}{4\pi b_n}\right)^m=C'_{d,\lambda} n^{\frac{k+1}{d}},
\end{equation}
where $C'_{d,\lambda}$ is also a constant only depending on $d$ and $\lambda$.

\vspace{1em}

 Let us recall Jackson's approximation
theorem (see  Appendix \ref{appendix}).

Let $f(x)$ be a $C^k$ and  $\Z^d$-periodic function on $\R^d$; then, for every
$N\in\N$ there exists a trigonometric polynomial $p_N(x)$ of degree at most
$N$ such that
\begin{equation}\label{jca}
\|p_N(x)-f(x)\|_{C^0}\leq C'_d2^kN^{-k}||f(x)||_{C^k},
\end{equation}
where  $C'_d$ is an absolute constant only depending on $d$ (independent of $f$, $k$ and $N$). Moreover, if $\int_{[0,1]^d}f(x)dx=0$, then $\int_{[0,1]^d}p_N(x)dx=0$.

Fix $n\in \N$. We need to construct $p_N(x)$ such that
\begin{equation}\label{need}
\|p_N(x)-f(x)\|_{C^0}\leq \frac{1}{4n},
\end{equation}
which implies
\begin{equation}\label{pn}
\max p_N(x)\leq \frac{1}{n},\quad -\min p_N(x)\in [\Delta,\Delta+\frac{1}{2n}].
\end{equation}
To achieve
\begin{equation}\label{need1}
\|p^n_N(x)-f_n(x)\|_{C^0}\leq \frac{1}{4n},
\end{equation}
 it suffices to require
\[C'_d2^kN^{-k}||f_n(x)||_{C^k}\leq \frac{1}{4n},\]
which together with (\ref{fck}) yields
\begin{equation}\label{N}
N\geq 2n^{\frac{1}{d}+\frac{1+d}{dk}}\cdot (4C'_dC'_{d,\lambda})^{\frac{1}{k}}.
\end{equation}
Note that $f_n\in C^\infty$. Take $k$ large enough, then (\ref{N}) can be verified by choosing
\[N=\left\lfloor n^{\frac{1}{d}+\eps}\right\rfloor.\]

Finally, we choose \[\phi_n(x)=\frac{-\Delta}{\min_{x\in[0,1]^d}p_N^n(x)}p_N^n(x),\quad x\in \R^d.\]
A direct calculation shows $\|\phi_n\|_{C^0}=O(\frac{1}{n})$ and $\|\phi_n\|_{C^1}=O(1)$. By the interpolation inequality (see \cite[Lemma 5]{Sa}), we have $\|\phi_n\|_{C^{1-\eps}}=O(\frac{1}{n^{\eps}})$.

This completes the proof of Lemma \ref{Csp}.
\End

In order to complete the proof of Theorem \ref{Mt}, we take

\begin{align}\label{clamb}
\Delta(\lambda):=\left\{
        \begin{array}{ll}
        1, & \lambda\in (0,\frac{1}{2}],\\
         & \\
        \frac{4}{3}(1-\lambda^2), & \lambda\in (\frac{1}{2},1).
         \end{array}
         \right.
\end{align}
It is clear to see that $\Delta(\lambda)\leq 1$ for all $\lambda\in (0,1)$, and $\Delta(\lambda)=O(1-\lambda)$ as $\lambda\to 1^-$.

Let $D\varphi_n^\lambda(x)=\phi_n(x)$ in Lemma \ref{Csp} with $\Delta=\Delta(\lambda)$. It follows that $\varphi^\lambda_n(x)=\int_{0}^x\phi_ndx +C_n$.  By choosing  suitable constants $C_n$, we have $\int_0^1\phi^\lambda_n(x)dx=0$. A direct calculation shows all of itema (I)-(III) are valid. Moreover, the criterion ($\diamondsuit$) is verified.

\vspace{1em}

\section{Proof of Theorem \ref{Mt2}}\label{S2}

Now we want to discuss a higher dimensional version of Theorem \ref{Mt} in the case of conformal symplectic maps of $\T^d\times \R^d$, as introduced in Definition \ref{confsympmap}.

Let us discuss some preliminary properties.

 Let $f$: $\T^d\times\R^d\mapsto \T^d\times\R^d$ be a $C^1$ exact conformal symplectic
 twist map (see Definition \ref{confsympmap}); in particular, let $F$ be its lift to $\R^{d}\times \R^d$ with  generating function $S:\R^{d}\times \R^d\to\R$. Then, $S$ is of class $C^2$ and $F$
is generated by the following equations:
\begin{equation}\label{geg}
\begin{cases}
\lambda y=-\frac{\pa S}{\pa x}(x,X),\\
Y=\frac{\pa S}{\pa X}(x,X),
\end{cases}
\end{equation}
where $F(x,y)=(X,Y)$. Assume that $S$ satisfies these extra conditions:
\begin{equation}\tag{H1}\label{H1}
\frac{\partial^2 S}{\partial x^2}(x,X)>0,\quad \frac{\partial^2 S}{\partial X^2}(x,X)>0, \quad \frac{\partial^2 S}{\partial x\partial X}(x,X)<0
\qquad \; (x,X)\in \R^{d}\times \R^d
\end{equation}
and
\begin{equation}\tag{H2}\label{H2}
\lim_{\|x-X\|\to+\infty}\frac{S(x,X)}{\|x-X\|}=+\infty.
\end{equation}
The following result is well known under these conditions (H1) and (H2). For $\lambda=1$, it was proved by Herman in \cite[Theorem 8.14]{H2} (see also \cite[Proposition A.1]{ams}).
\begin{Proposition}\label{diff}
\
\begin{itemize}
\item
Given $X\in \R^d$, the map $x\mapsto \frac{\partial S}{\partial X}(x,X)$ is a diffeomorphism on $\R^d$.
\item
Given $x\in \R^d$, the map $X\mapsto \frac{\partial S}{\partial x}(x,X)$ is a diffeomorphism on $\R^d$.
\end{itemize}
\end{Proposition}

\medskip

Let us now discuss some \textit{a-priori} Lipschitz estimates for invariant graphs.

\begin{Proposition}\label{lipest}
Let $F:\R^{2d}\to\R^{2d}$ be a conformal symplectic
 twist map generated by $S$ satisfying (\ref{H1}) and (\ref{H2}). If $F$ admits a $C^0$-invariant Lagrangian graph $\mathcal{L}=\{(x,\tilde{\Psi}(x))\ |\ x\in \R^{d}\}$, then
$\tilde{\Psi}$ is Lipschitz continuous, and satisfies
\begin{equation}\label{lipineq}
\|D\tilde{\Psi}\|_{L^\infty}\leq \sup_{x\in \R^d}\left\{\frac{1}{\lambda}\|\frac{\partial^2 S}{\partial x^2}(x,\tilde{\Psi}(x))\|_{\infty}, \|\frac{\partial^2 S}{\partial X^2}(x,\tilde{\Psi}(x))\|_{\infty}\right\},
\end{equation}
where $\|D\tilde{\Psi}\|_{L^\infty}:=\mathrm{ess\ sup}_{x\in\R^d}\|D\tilde{\Psi}(x)\|_{\infty}$.
\end{Proposition}

\Proof Let $\tilde{\eta}:\R^d\to\R$ be the lift of $\eta:\T^d\to\R$ in Definition \ref{lag}. Then for all $x\in \R^d$,
\[ \til{\Psi}(x)=c+D\til{\eta}(x).\]
Let $\hat{\Psi}(x):=\langle c,x\rangle+\til{\eta}$. Then $D\hat{\Psi}(x)=\til{\Psi}$.
Denote
\[K(x,X):=S(x,X)+\lambda\hat{\Psi}(x)-\hat{\Psi}(X).\]
If $F$ admits a $C^0$-invariant Lagrangian graph $\mathcal{L}=\{(x,\tilde{\Psi}(x))\ |\ x\in \R^{d}\}$, then there exists a homeomorphism $g:\R^d\to\R^d$ such that
\[F(x,\til{\Psi}(x))=(g(x),\til{\Psi}(g(x))),\quad F(g^{-1}(x),\til{\Psi}(g^{-1}(x)))=(x,\til{\Psi}(x)),\]
which together with (\ref{geg}) implies
\begin{equation}\label{equal}
\lambda\til{\Psi}(x)=-\frac{\pa S}{\pa x}(x,g(x))=\lambda \frac{\pa S}{\pa X}(g^{-1}(x),x).
\end{equation}
We claim the following statements.
\begin{itemize}
\item [($\dag$)] Given $X\in \R^d$, the map $x\mapsto K(x,X)$ has a unique minimal point $x_0=g^{-1}(X)$.
\item [($\ddag$)] Given $x\in \R^d$, the map $X\mapsto K(x,X)$ has a unique minimal point $X_0=g(x)$.
\end{itemize}
We only prove item (1), since item (2) follows from a similar argument.
Given $X\in \R^d$, let $x_0$ be a critical point of the map $x\mapsto K(x,X)$. Then
\[\frac{\pa K}{\pa x}(x_0,X)=\frac{\pa S}{\pa x}(x_0,X)+\lambda\til{\Psi}(x_0)=0.\]
By (\ref{equal}),
\[\frac{\pa K}{\pa x}(x_0,g(x_0))=\frac{\pa S}{\pa x}(x_0,g(x_0))+\lambda\til{\Psi}(x_0)=0.\]
 It follows from Proposition \ref{diff} that $g(x_0)=X$. Then $x_0=g^{-1}(X)$, since $g:\R^d\to\R^d$ is a homeomorphism. By (H2),
 \begin{itemize}
\item  given $X\in \R^d$, $K(x,X)\to +\infty$ as $\|x\|\to+\infty$;
\item given $x\in \R^d$, $K(x,X)\to +\infty$ as $\|X\|\to+\infty$.
\end{itemize}
 Then  $x_0=g^{-1}(X)$ is the unique minimal point, where  uniqueness follows from the fact that $K(x,X)$ is strictly convex for a given $X$
as it follows from the assumption $\frac{\partial^2 S}{\partial x^2}(x,X)>0$ in \eqref{H1}.
   This completes the proof of the claim.

Given $s\in \R$, $v\in \R^d$, $x,X\in \R^d$, denote
\[\hat{\Delta}_{sv}(x):=\frac{1}{s^2}\left(K(x+sv,g(x))+K(x-sv,g(x))-2K(x,g(x))\right),\]
\[\check{\Delta}_{sv}(X):=\frac{1}{s^2}\left(K(g^{-1}(X),X+sv)+K(g^{-1}(X),x-sv)-2K(g^{-1}(X),X)\right).\]
Based on the claims ($\dag$) and ($\ddag$), we have $\hat{\Delta}_{sv}\geq 0$, and $\check{\Delta}_{sv}\geq 0$. By the definition of $K$,
\begin{align*}
\hat{\Delta}_{sv}(x):=&\frac{1}{s^2}\left(S(x+sv,g(x))+S(x-sv,g(x))-2S(x,g(x))\right)\\
&+\frac{\lambda}{s^2}\left(\hat{\Psi}(x+sv)+\hat{\Psi}(x-sv)-2\hat{\Psi}(x)\right).
\end{align*}
Let $s\to 0$. In the sense of distribution, we have
\begin{equation}\label{dis1}
\lambda v^TD^2\hat{\Psi}(\cdot)v+v^T\frac{\pa^2 S}{\pa x^2}(\cdot,g(\cdot))v=:\nu^1_v\geq 0,
\end{equation}
where $v^T$ denotes the transpose of $v\in \R^d$. Similarly, we have
\begin{equation}\label{dis2}
v^T\frac{\pa^2 S}{\pa X^2}(g^{-1}(\cdot),\cdot)v- v^TD^2\hat{\Psi}(\cdot)v=:\nu^2_v\geq 0.
\end{equation}
By the same argument as \cite[Page 64-65]{H2}, we know that for $|v|>0$, $\nu^1_v$ and $\nu^2_v$ are positive Radon measures, and $x\mapsto v^TD\til{\Psi}(x)v\in L^{\infty}$. It follows that $\til{\Psi}$ is Lipschitz continuous on $\R^d$. Then $\til{\Psi}$ is differentiable on a full Lebesgue measure set $\mathfrak{D}$. By (\ref{dis1}) and (\ref{dis2}), for each $x,X\in \mathfrak{D}$,
\begin{equation}\label{dis3}
\lambda D\til{\Psi}(x)+\frac{\pa^2 S}{\pa x^2}(x,g(x))\geq 0,\quad \frac{\pa^2 S}{\pa X^2}(g^{-1}(X),X)- D\til{\Psi}(X)\geq 0.
\end{equation}
By Proposition \ref{diff}, it follows from (\ref{equal}) that $g$ is bi-Lipschitz continuous. Then $\bar{\mathfrak{D}}:=\mathfrak{D}\cap g(\mathfrak{D})$ is also a full Lebesgue measure set in $\R^d$. By (\ref{dis3}), for each $x\in \bar{\mathfrak{D}}$,
\begin{equation}\label{dis3+}
\lambda D\til{\Psi}(x)+\frac{\pa^2 S}{\pa x^2}(x,g(x))\geq 0,\quad \frac{\pa^2 S}{\pa X^2}(g^{-1}(x),x)- D\til{\Psi}(x)\geq 0.
\end{equation}
This yields (\ref{lipineq}).
\End

\subsection{Herman's formula for conformal symplectic maps on $\T^d$} \label{subsechermantd}
We can now prove an analogue of Herman's formula proved in Section \ref{Herman1} in this setting.
By (\ref{fbeta}), $F_{\lambda,\beta}$ is generated by
\begin{equation}\label{S0}
S_{\lambda,\beta}(x,X):=\frac{1}{2}\langle X-x,X-x\rangle-\lambda \langle \beta,X-x\rangle,
\end{equation}
where $\langle \cdot,\cdot\rangle$ denotes the standard inner product in $\R^d$.

Let $\Phi$ be a $\Z^d$-periodic $C^2$ function  satisfying $\int_{[0,1]^d}\Phi(x)dx=0$. The notation $[0,1]^d$ denotes the $d$-fold Cartesian product of $[0,1]$. Let
\[S^\flat_{\lambda,\beta}(x,X):=\frac{1}{2}\langle X-x,X-x\rangle-\lambda \langle \beta,X-x\rangle+\lambda \Phi(x).\]
Let $F_{\lambda,\beta}^\flat:\R^{2d}\to\R^{2d}$ be generated by $S^\flat_{\lambda,\beta}$.
It follows that
\[F_{\lambda,\beta}^\flat(x,y)=(x+\lambda (\beta+y+D\Phi(x)),\lambda (y+D\Phi(x))).\]

Similar to Proposition \ref{hf1}, we have
\begin{Proposition}\label{hf2}
$F^\flat_{\lambda,\beta}$ admits a $C^0$-invariant Lagrangian graph $\tilde{\Gamma}:=\{(x,\tilde{\Psi}(x))\ |\ x\in \R^d\}$ if and only if the following holds
\begin{equation}\tag{B}\label{bb}
\frac{1}{1+\lambda}g(x)+\frac{\lambda}{1+\lambda}g^{-1}(x)=x+\frac{\lambda}{1+\lambda}((1-\lambda)\beta+D\Phi(x)) \qquad \forall\; x\in \R^d,
\end{equation}
where $g(x):=x+\lambda(\beta+\tilde{\Psi}(x)+D\Phi(x))$.
\end{Proposition}

{The proof is similar to the one of Proposition \ref{hf1} and we omit it.}\\

\begin{Remark}
For $\lambda=1$, the formula (B) was established by Herman in 1990 (see \cite{H3}).
\end{Remark}

{
\begin{Remark}\label{RemarkmatrixA2}
In the setting introduced in Remark \ref{RemarkmatrixA1}, namely for maps of the form
\[F_{\lambda,\beta}^\flat(x,y)=\left(x+\lambda (\beta+D\Phi(x)+A^{-1}y),\lambda (y+AD\Phi(x))\right),\]
where $A$ is a $d\times d$ symmetric positive definite  matrix,  $\beta\in \R^d$ be a constant vector, a direct calculation shows that
\[\left(F_{\lambda,\beta}^\flat\right)^{-1}(x,y)=\left(x-\lambda\beta-A^{-1}y,\frac{1}{\lambda}y-D\Phi(x-\lambda\beta-A^{-1}y)\right).\]
Condition \eqref{bb} in Proposition \ref{hf2} remains the same, substituting $g(x):=x+\lambda(\beta+D\Phi(x)+ A^{-1}\tilde{\Psi}(x))$.
The remaining argument for proving Theorem \ref{Mt2} in this setting is exactly the same as what we are going to show hereafter.\\
\end{Remark}
}

By Proposition \ref{lipest}, $\tilde{\Psi}$ is Lipschitz continuous. Then $g$ is bi-Lipschitz {(proved in the proof of Proposition \ref{lipest})}. Let $\mathfrak{D}$ be the set of points of differentiability of $g$,
then $\mathfrak{D}$ has full Lebesgue measure on $\R^d$. By the construction of $F^\flat_{\lambda,\beta}$, for each $x\in \mathfrak{D}$, $Dg(x)$ and $Dg^{-1}(x)$ are symmetric positive definite matrices. We differentiate
(\ref{bb}) on $\mathfrak{D}$ to get
\begin{equation}
\frac{1}{1+\lambda}Dg(x)+\frac{\lambda}{1+\lambda}(Dg)^{-1}(g^{-1}(x))=\mathbb{I}_d+\frac{\lambda}{1+\lambda}D^2\Phi(x),
\end{equation}
where $\mathbb{I}_d$ is the $d\times d$ identity matrix. Let $G(x):=Dg(x)$. Then $G^{-1}(g^{-1}(x))=Dg^{-1}(x)$. Denote $E(x):=D^2\Phi(x)$. It follows that
\begin{equation}\label{deri-b}
\frac{1}{1+\lambda}\cdot\frac{1}{d}\mathrm{tr}(G(x))+\frac{\lambda}{1+\lambda}\cdot\frac{1}{d}\mathrm{tr}(G^{-1}(g^{-1}(x)))=1+\frac{\lambda}{1+\lambda}\cdot
\frac{1}{d}\mathrm{tr}(E(x)).
\end{equation}
Let us first introduce some notation:
\begin{itemize}
\item Denote $T(x):=\frac{1}{d}\mathrm{tr}(E(x))$,
\[M:=\max_{x\in [0,1]^d}T(x),\quad m:=-\min_{x\in [0,1]^d}T(x).\]
\item Denote
$\mathfrak{D}_0:=\mathfrak{D}\cap [0,1]^d$,
\[G_+:=\max\left\{\sup_{\mathfrak{D}_0}\frac{1}{d}\mathrm{tr}(G(x)),\sup_{\mathfrak{D}_0}\frac{1}{d}\mathrm{tr}(G^{-1}(x))\right\},\]
\[G_-:=\min\left\{\inf_{\mathfrak{D}_0}\frac{1}{d}\mathrm{tr}(G(x)),\inf_{\mathfrak{D}_0}\frac{1}{d}\mathrm{tr}(G^{-1}(x))\right\}.\]
\end{itemize}

\smallskip

Note that for each $x\in \mathfrak{D}$, $Dg(x)$ and $Dg^{-1}(x)$ are symmetric positive definite matrices. In fact,  by definition, $g(x):=x+\lambda(\beta+\tilde{\Psi}(x)+D\Phi(x))$. Then $Dg(x)=\mathbb{I}_d+\lambda D\tilde{\Psi}+D^2\Phi(x)$ and $\tilde{\Psi}=c+D\tilde{\eta}$.
\begin{itemize}
\item {\it Symmetry:}
We only need to prove that $D^2\tilde{\eta}$ is symmetric. In fact, by Proposition \ref{lipest}, $\tilde{\eta}$ is of class $C^{1,1}$. It follows that $D^2\tilde{\eta}(x)$ is symmetric for each $x\in \mathfrak{D}$.
\item {\it Positive definiteness:} By Proposition \ref{lipest}, $g$ is bi-Lipschitz continuous. Then $Dg(x)$ and $Dg^{-1}(x)$ are non-degenerate for each $x\in \mathfrak{D}$. According to (\ref{equal}), we have
\[\lambda D\til{\Psi}(x)+\frac{\pa^2 S}{\pa x^2}(x,g(x))=-\frac{\partial^2 S}{\partial x\partial X}(x,g(x))Dg(x),\]
\[\frac{\pa^2 S}{\pa X^2}(g^{-1}(x),x)- D\til{\Psi}(x)=-\frac{\partial^2 S}{\partial x\partial X}(g^{-1}(x),x)Dg^{-1}(x).\]
    Recalling (H1), we assume
 $\frac{\partial^2 S}{\partial x\partial X}(x,X)<0$. It follows from  (\ref{dis3+}) that  $Dg(x)$ and $Dg^{-1}(x)$ are positive definite for  each $x\in \mathfrak{D}$.\\
\end{itemize}

\medskip

Then $0<G_-\leq G_+$. Let us apply the  AM-GM inequality. Note that
 \[\mathrm{tr}(G(x))=\lambda_1+\cdots+\lambda_d,\qquad {\rm and} \qquad \mathrm{tr}(G^{-1}(x))=\frac{1}{\lambda_1}+\cdots+\frac{1}{\lambda_d}.\]
 Then, we have
\[\mathrm{tr}(G(x))\mathrm{tr}(G^{-1}(x))\geq d^2 \qquad {\rm and} \qquad
G_-\geq \frac{1}{G_+},\]
which imply
\begin{equation}\label{ggGG2}
\frac{1}{G_+}\leq 1+\frac{\lambda}{1+\lambda}m,
\end{equation}
which still implies $m>-1-\lambda$.
The following part is also divided into two cases.
\begin{itemize}
\item {\bf{Case 1:}} $G_+=\sup_{\mathfrak{D}_0}\frac{1}{d}\mathrm{tr}(G(x))$;
\item {\bf{Case 2:}} $G_+=\sup_{\mathfrak{D}_0}\frac{1}{d}\mathrm{tr}(G^{-1}(x))$.
\end{itemize}

Similar to the argument from (\ref{G}) to (\ref{mmMM2}), we have if Case 1 holds, $m>-1-\lambda$ and $M\to 0^+$, then
\[-m\leq \frac{1+\lambda}{1-\lambda}M+O(M^2).\]
On the other hand, if
we have if Case 1 holds, $m>-1-\lambda$ and $M\to 0^+$, then
\[-m\leq 1-\lambda^2+\frac{\lambda(1+\lambda)}{1-\lambda}M+O(M^2).\]

\medskip

\subsection{Construction of the perturbation}
Let us recall
\begin{align*}
\Delta(\lambda):=\left\{
        \begin{array}{ll}
        1, & \lambda\in (0,\frac{1}{2}],\\
         & \\
        \frac{4}{3}(1-\lambda^2), & \lambda\in (\frac{1}{2},1).
         \end{array}
         \right.
\end{align*}

Let $T_n(x)=\phi_n(x)$ in Lemma \ref{Csp} with $\Delta=\Delta(\lambda)$.
Note that $\Phi_n$ as the perturbation of the generating function $S_{\lambda,\beta}$ (see (\ref{S0})) is uniquely determined by
\[\frac{1}{d}\triangle \Phi_n(x)=T_n(x).\]
Then the proof of Theorem \ref{Mt2} can be completed by a similar argument as in Theorem \ref{Mt}.

 \vspace{2em}

\setcounter{equation}{0}

\renewcommand\theequation{A.\arabic{equation}} 
\appendix

\section{On multivariate Jackson's  approximation theorem}\label{appendix}

In this section we present a proof of Jackson's approxation theorem, using an argument similar to  \cite[Theorem 2.12]{A}.\\

\begin{Remark}
Jackson's approximation has been also employed in \cite[Lemma 3]{Sa}, where it is presented in a slightly modified form. Specifically, the formulation in \cite{Sa} involves the $C^r$ norm of the difference on the left-hand side, whereas our analysis only requires the $C^0$ norm for that particular term. Furthermore, it appears that the significance of the degree $N$ in the trigonometric polynomial has not been sufficiently highlighted. Another crucial aspect is the dependence of the constant $C$ in Theorem \ref{appp}  solely on the dimension $d$. As noted in \cite{Z}, for the case $d=1$, the explicit computation of this constant $C$ was first achieved by Favard \cite{Fav}. \end{Remark}

\bigskip

 First of all, we need
some notations. Define
\[C^\infty_{\Z^d}(\R^d,\R):=\left\{f:\R^d\rightarrow\R|f\in C^\infty(\R^d,\R)\ \text{and}\ \Z^d-\text{periodic}\right\}.\]
 Let $f(x)\in C^\infty_{1}(\R^d,\R)$. The $m$-th
 Fej\'{e}r-polynomial of $f$ with respect to $x_j$ is given by
 \begin{equation}\label{B1}
F_m^{[j]}(f)(x):=\frac{2}{m\pi}\int_{-1/4}^{1/4}f(x+2te_j)\left(\frac{\sin(2\pi mt)}{\sin (2\pi
t)}\right)^2dt,
 \end{equation}
where $x\in\R^d$, $m\in\N$, $j\in\{1,\ldots,d\}$ and $e_j$ is the
$j$-th vector of the canonical basis of $\R^d$. $F_m^{[j]}(f)(x)$ is
a trigonometric polynomial in $x_j$ of degree at most $m-1$. By
\cite{Z},
\[\frac{2}{m\pi}\int_{-1/4}^{1/4}\left(\frac{\sin(2\pi mt)}{\sin (2\pi
t)}\right)^2dt=1,\] hence, from (\ref{B1}), we have
\[||F_m^{[j]}(f)||_{C^0}\leq ||f||_{C^0}.\]
We denote
\[P_m^{[j]}(f):=2F_{2m}^{[j]}(f)-F_m^{[j]}(f).\] It is easy to see
that $P_m^{[j]}(f)$ is a trigonometric polynomial in $x_j$ of degree
at most $2m-1$. Moreover,
 \begin{equation}\label{B2}
||P_m^{[j]}(f)||_{C^0}\leq 3||f||_{C^0},
 \end{equation}
 \begin{equation}\label{B3}
P_m^{[j]}(af+bg)=aP_m^{[j]}(f)+bP_m^{[j]}(g),
 \end{equation}
where $a,b\in\R$ and $f,g\in C^\infty_{1}(\R^d,\R)$. For
$k\in\{1,\ldots,d\}$, $j_1,\ldots,j_k\in\{1,\ldots,d\}$ with
$j_p\neq j_q$ for $p\neq q$. Let $m_1,\ldots,m_k\in\N$ and $f\in
C^\infty_{1}(\R^d,\R)$, we define
 \begin{equation}\label{B4}
P_{m_1,\ldots,m_k}^{[j_1,\ldots,j_k]}(f):=P_{m_1}^{[j_1]}\left(P_{m_2}^{[j_2]}\left(\cdots\left(P_{m_k}^{[j_k]}(f)\right)\cdots\right)\right).
 \end{equation}
It is easy to see that for all $l\in \{1,\ldots,k\}$,
$P_{m_1,\ldots,m_k}^{[j_1,\ldots,j_k]}(f)$ are trigonometric
polynomials in $x_{j_l}$ of degree at most $2m_l-1$, also known as generalized de la Vall\'{e}e Poussin polynomial.
\begin{Theorem}\label{appp}
Let $f\in C^\infty_{1}(\R^d,\R)$, $r_1,\ldots,r_d\in\N$,
$m_1,\ldots,m_d\in\N$, then we have
\begin{equation}\label{B8}
\|f-P_{m_1,\ldots,m_d}^{[1,\ldots,d]}(f)\|_{C^0}\leq
C\sum_{j=1}^{d}\frac{1}{{m_j}^{r_j}}\left\|\frac{\partial
^{r_j}f}{\partial {x_j}^{r_j}}\right\|_{C^0},
\end{equation}
where $C$ is a constant only depending on $d$.
\end{Theorem}

\Proof We will prove Theorem \ref{appp} by induction.
The case $d=1$ is covered by the classical Jackson's approximation theorem after Favard \cite{Fav}. See \cite[Theorem 13.6, p115 and Notes, p377]{Z}). More
precisely, for $f\in C^\infty_{1}(\R,\R)$, $m,r\in\N$, we have
\begin{equation}\label{B5}
\|f-P_m^{[1]}(f)\|_{C^0}\leq
C_1\frac{1}{{m}^{r}}\left\|\frac{\partial ^{r}f}{\partial
{x}^{r}}\right\|_{C^0},
\end{equation}
where $C_1$ is an absolute constant independent of $f$, $m$ and $r$. Let the assertion be true for $d=k\in\N$. We verify it for $d=k+1$.
Consider the functions $f(x_1,\cdot)$ with $x_1$ as a real
parameter. Then by the assertion for $d$, we have
\[\|f(x_1,\cdot)-P_{m_2,\ldots,m_{k+1}}^{[2,\ldots,k+1]}(f)(x_1,\cdot)\|_{C^0}\leq
C_k\sum_{j=2}^{k+1}\frac{1}{{m_j}^{r_j}}\left\|\frac{\partial
^{r_j}f}{\partial {x_j}^{r_j}}\right\|_{C^0},\] hence,
\begin{equation}\label{B6}
\|f-P_{m_2,\ldots,m_{k+1}}^{[2,\ldots,k+1]}(f)\|_{C^0}\leq
C_k\sum_{j=2}^{k+1}\frac{1}{{m_j}^{r_j}}\left\|\frac{\partial
^{r_j}f}{\partial {x_j}^{r_j}}\right\|_{C^0}.
\end{equation}
Let $\hat{x}_j\in\R^k$ denote the vector $x\in\R^{k+1}$ without its
$j$-th entry. For the functions $f(\cdot,\hat{x}_1)$, from
(\ref{B5}), it follows that
\[\|f(\cdot,\hat{x}_1)-P_{m_1}^{[1]}(f)(\cdot,\hat{x}_1)\|_{C^0}\leq C_1\frac{1}{{m_1}^{r_1}}\left\|\frac{\partial ^{r_1}f}{\partial
{x_1}^{r_1}}\right\|_{C^0},\]hence,
\begin{equation}\label{B7}
\|f-P_{m_1}^{[1]}(f)\|_{C^0}\leq
C_1\frac{1}{{m_1}^{r_1}}\left\|\frac{\partial ^{r_1}f}{\partial
{x_1}^{r_1}}\right\|_{C^0}.
\end{equation}
By (\ref{B2}), (\ref{B3})),(\ref{B4}) and (\ref{B6}), we have
\begin{align*}
\left\|P_{m_1}^{[1]}(f)-P_{m_1,\ldots,m_{k+1}}^{[1,\ldots,k+1]}(f)\right\|_{C^0}&=\left\|P_{m_1}^{[1]}(f)-P_{m_1}^{[1]}\left(P_{m_2,\ldots,m_{k+1}}^{
[2,\ldots,k+1]}(f)\right)\right\|_{C^0},\\
&=\left\|P_{m_1}^{[1]}\left(f-P_{m_1}^{[1]}P_{m_2,\ldots,m_{k+1}}^{
[2,\ldots,k+1]}(f)\right)\right\|_{C^0},\\
&\leq 3\left\|f-P_{m_1}^{[1]}P_{m_2,\ldots,m_{k+1}}^{
[2,\ldots,k+1]}(f)\right\|_{C^0},\\
&\leq 3C_k\sum_{j=2}^{k+1}\frac{1}{{m_j}^{r_j}}\left\|\frac{\partial
^{r_j}f}{\partial {x_j}^{r_j}}\right\|_{C^0},
\end{align*}
which together with (\ref{B7}) implies that
\begin{align*}
\left\|f-P_{m_1,\ldots,m_{k+1}}^{[1,\ldots,d+1]}(f)\right\|_{C^0}&\leq \left\|f-P_{m_1}^{[1]}(f)\right\|_{C^0}+\left\|P_{m_1}^{[1]}(f)-P_{m_1,\ldots,m_{d+1}}^{[1,\ldots,k+1]}(f)\right\|_{C^0},\\
&=\left\|f-P_{m_1}^{[1]}(f)\right\|_{C^0}+\left\|P_{m_1}^{[1]}(f)-P_{m_1}^{[1]}\left(P_{m_2,\ldots,m_{k+1}}^{
[2,\ldots,k+1]}(f)\right)\right\|_{C^0},\\
&\leq C_1\frac{1}{{m_1}^{r_1}}\left\|\frac{\partial
^{r_1}f}{\partial
{x_1}^{r_1}}\right\|_{C^0}+3C_k\sum_{j=2}^{k+1}\frac{1}{{m_j}^{r_j}}\left\|\frac{\partial
^{r_j}f}{\partial {x_j}^{r_j}}\right\|_{C^0},\\
&\leq
C_{k+1}\sum_{j=1}^{k+1}\frac{1}{{m_j}^{r_j}}\left\|\frac{\partial
^{r_j}f}{\partial {x_j}^{r_j}}\right\|_{C^0}.
\end{align*}
This finishes the proof of Theorem \ref{appp}. \End

We choose $m_1=\ldots=m_d=\bar{m}$. Let $\bar{r}$ be the value  such that
\[\frac{1}{\bar{m}^{\bar{r}}}\left\|\frac{\partial
^{r_{\bar{j}}}f}{\partial
{x_{\bar{j}}}^{r_{\bar{j}}}}\right\|_{C^0}=\max_{1\leq j \leq
d}\left\{\frac{1}{{\bar{m}}^{r_j}}\left\|\frac{\partial
^{r_j}f}{\partial {x_j}^{r_j}}\right\|_{C^0}\right\}.\] Hence, we
have
\begin{align*}
\|f-P_{m_1,\ldots,m_d}^{[1,\ldots,d]}(f)\|_{C^0}&\leq
dC_d\frac{1}{\bar{m}^{\bar{r}}}\left\|\frac{\partial
^{r_{\bar{j}}}f}{\partial
{x_{\bar{j}}}^{r_{\bar{j}}}}\right\|_{C^0},\\
&\leq
C'_d\frac{1}{\bar{m}^{\bar{r}}}\|f\|_{C^{\bar{r}}}.
\end{align*}
For the simplicity of notations, we denote
\[p_N(x)=P_{\bar{m},\ldots,\bar{m}}^{[1,\ldots,d]}(f)(x),\]
where $x=(x_1,\ldots,x_d)$ and $N=2\bar{m}-1$. Moreover, we
denote $k:=\bar{r}$, then
\begin{equation}\label{B9}
\|f(x)-p_N(x)\|_{C^0}\leq C2^kN^{-k}\|f(x)\|_{C^k},
\end{equation}
where $C$ is a constant only depending on $d$.\\

\end{document}